\documentclass{amsart}
\usepackage{graphicx}
\theoremstyle{definition}

\theoremstyle{remark}

\numberwithin{equation}{section}

\begin{document}

\title[]{Identities for hyperelliptic $\wp$-functions of genus one, two and three in covariant form.}
\author{Chris Athorne}%
\address{Maths. Glasgow}%
\email{ca@maths.gla.ac.uk}%


\begin{abstract}
We give a covariant treatment of the quadratic differential
identities satisfied by the $\wp$-functions on the Jacobian of
smooth hyperelliptic curves of genus $\leq 3$.
\end{abstract}
\maketitle
\section{Introduction}
A classical problem in the theory of a planar $(n,s)$ algebraic
curve is a description of the differential equations satisfied by
meromorphic, multiply periodic functions defined on its Jacobian
variety. In the genus $g$ hyperelliptic case ($n=2,\,s=2g+2$) the
field of such functions is entirely described in terms of certain
$\wp_{ij}$ functions which generalize the Weierstrass $\wp$-function
on the elliptic curve, the genus one case.

The derivation of these identities has been a major concern over the
last ten to fifteen years and many results have been published: see
\cite{BEL1997,EEL2000,EEP2003} for seminal literature.

The aim of this paper is to promote a new methodology which
considerably simplifies the derivation and presentation of these
identities by utilizing elementary representation theory. The
fundamental observation is that the underlying algebraic curves
belong to generic families permuted under an ${\mathfrak sl}_2$
action. This can be interpreted \cite{AEE2003,AEE2004} as a
covariance property that translates into covariance of the
$\wp$-function identities. This means that each polynomial identity
between derivatives of the $\wp$-function belongs to a finite
dimensional representation of ${\mathfrak sl}_2$, the knowledge of
which depends only upon a highest weight element. It is only
necessary to find these highest weight identities to generate the
other identities in the representation.

However, a requirement of this approach is that we develop the
theory for the generic member of the family of curves. This is in
contrast to former treatments where a simpler, normal form is
exploited by moving a branch point to infinity, i.e. removing the
highest degree term.

The only case where the covariant equations are written down is for
genus two hyperelliptic curves by Baker \cite{B1897}. He achieves
this by establishing the equations for the curve in normal form and
then undoing the ``normalizing" transformation's effect on the
identities. Even so he finds it necessary to introduce a ``fudge
factor" to restore full covariance.

This ``fudge factor" points to another problem. Not only must the
curve be in general position but the fundamental (Kleinian)
definition of the $\wp$-function \cite{B1897,B1903,B1907,BEL1997}
must itself be rendered covariant. This problem reasserts itself in
the next highest genus and the Baker equations for the genus three
curve \cite{B1903} are nowhere written down in covariant form.

The same issues occur in purely algebraic treatments, that of
Cassels and Flynn for instance \cite{CF1996}. The formulation of
their approach, important for curves over general fields, can also
be rendered covariant and will be discussed in another publication.
In this paper we work entirely over $\mathbb C$.

In this respect a note on the approach of the papers
\cite{AEE2003,AEE2004} by the present author and collaborators is in
order. What was attempted in those papers was a radically different
approach to the analytic theory based on a very simple definition of
the $\wp$-function, quite different to Klein's but with some
philosophical proximity to that of \cite{CF1996}. However, whilst
this was an effective approach to genus two, attempts so far to
extend it to higher genus have foundered on finding the
corresponding, simple definition of the $\wp$-function.

The programme of the current paper is, therefore, firstly to define
the $\wp$ function in a covariant way and secondly to derive the
identities it satisfies by combining the traditional technique of
expansion about a chosen point with the Lie algebraic representation
theory. We do this for genera one, two and three to recover known
sets of differential equation or their equivalents. The emphasis is
placed on the methodology.

The results so obtained are rather beautiful generalizations of the
formulae found in \cite{B1903,B1907,BEL1997}. Most of all we obtain
a covariant bordered determinantal form of the set of quadratic
identities in the $\wp_{ijk}$ for genus two, familiar from Baker's
work \cite{B1907}, and a new generalization of this formula to the
genus three case involving a doubly bordered determinant. These
quadratic relations should presumably be regarded as the most
fundamental differential identities and it is a positive feature of
the covariant machinery that it produces them in a systematic manner
at the simplest level.

\section{Lie algebraic operations}

Curves of the form
\begin{equation}
v(x,y;a_0,\ldots,a_{2g+2)})=y^2-\sum_{i=0}^{2g+2}\binom{2g+2}{i}a_ix^i=0
\end{equation}
are generically hyperelliptic and of genus $g$: that is, unless some
special relations obtain between the coefficients.

The family of such curves is permuted under transformations given by
\begin{eqnarray}
x\mapsto X=\frac{\alpha x+\beta}{(\gamma x+\delta)},\\
y\mapsto Y=\frac{y}{(\gamma x +\delta)^{g+1}},
\end{eqnarray}
where $$\alpha\delta-\beta\gamma=1,$$ mapping the above curve into
\begin{equation}
V(X,Y;A_0,\ldots,A_{2g+2})=Y^2-\sum_{i=0}^{2g+2}\binom{2g+2}{i}A_iX^i=0
\end{equation}
the $A_i$ being functions of the $a_i$ and the parameters
$\alpha,\beta,\gamma$ and $\delta$.

This can be restated as infinitesimal \emph{covariance} conditions,
\begin{eqnarray}
{\bf e}v(x,y;a_0,\ldots,a_{2g+2})=0\\
{\bf
f}v(x,y;a_0,\ldots,a_{2g+2})+2(g+1)xv(x,y;a_0,\ldots,a_{2g+2})=0
\end{eqnarray}
where the generators $\bf e$ and $\bf f$ are given by
\begin{eqnarray}
{\bf e}=\partial_x-\sum_{i=0}^{2g+2}(2g+2-i)a_{i+1}\partial_{a_i}\\
{\bf
f}=-x^2\partial_x-(g+1)xy\partial_y-\sum_{i=0}^{2g+2}ia_{i-1}\partial_{a_i}\\
{\bf
h}=-2x\partial_x-(2g+2)y\partial_y-\sum_{i=0}^{g+1}ia_{i}\partial_{a_i}.
\end{eqnarray}

These generators satisfy the ${\mathfrak sl}_2$ commutation
relations,
\begin{equation}
[{\bf h},{\bf e}]=2{\bf e},\quad [{\bf h},{\bf f}]=-2{\bf f},\quad
[{\bf e},{\bf f}]={\bf h}.
\end{equation}

The coefficients $a_0,a_1,\ldots,a_{2g+2}$ are a basis for a $2g+3$
dimensional representation.

The space of holomorphic differentials on the curve is spanned by
the set \[\{\frac{x^{i-1}dx}{y}|i=1,\ldots,g\}\] and the symmetric
sums of each of these differentials taken over $g$ copies of the
curve,
\begin{equation}
du_i=\sum_{j=1}^g\frac{x_j^{i-1}dx_j}{y_j}
\end{equation}
are a basis for holomorphic one-forms on the Jacobian variety of the
curve.

One checks the following action of ${\mathfrak sl_2}$:
\begin{eqnarray}
{\bf e}du_i&=&(i-1)du_{i-1}\\
{\bf f}du_i&=&(g-i)du_{i+1}
\end{eqnarray}
and it then follows that
\begin{eqnarray}\label{ef}
{\bf e}\partial_{u_i}&=&-i\partial_{u_{i+1}}\\
{\bf f}\partial_{u_i}&=&-(g-i+1)\partial_{u_{i-1}}
\end{eqnarray}

\section{Covariant Klein relations}
Our starting point will be the Kleinian definition of the doubly
indexed $\wp$ functions: $\wp_{ij}=\wp_{ji}$ \cite{BEL1997}. The
indices are to be thought of as derivatives with respect to the
variables $u_i$. There are thus integrability conditions of the
form,
\begin{equation}
\wp_{ij,k}=\wp_{ik,j}=\wp_{kj,i}\quad\forall i,j,k.
\end{equation}
For the moment we think of these objects purely as indexed symbols
satisfying algebraic rules of differentiation and a set of
identities to be specified shortly. However they are not
traditionally defined in a covariant manner, that is in a way that
respects the further relations following from (\ref{ef}), namely,
\begin{eqnarray}
{\bf e}\wp_{ij}&=&-i\wp_{i+1\,j}-j\wp_{i\,j+1}\\
{\bf f}\wp_{ij}&=&-(g-i+1)\wp_{i-1\,j}-(g-j+1)\wp_{i\,j-1}
\end{eqnarray}

In order to proceed we need to adjust the fundamental definition by
adding correction terms without destroying the fundamental
singularity properties of the $\wp_{ij}$.

How to do this is best seen by example and we explain it now for the
genus two case.

The classical definition in genus two assumes a normal form with
branch point at infinity, $a_6=0, a_5=\frac23$, and is:
\begin{equation}
\wp_{11}+(x_i+x)\wp_{12}+xx_i\wp_{22}=\frac{F(x,x_i)-yy_i}{4(x-x_i)^2}
\end{equation}
where $i=1,2$ and $\wp$ is a function of the argument $\int^xd{\bf
u}+\int^{x_1}d{\bf u}+\int^{x_2}d{\bf u}$, ${\bf u}=(u_1,u_2)$. The
function $F(x,x_i)$ is the classical polar form
\begin{equation}
F(x,x_i)=2(x+x_i)x^2x^2_i+15a_4x^2x_i^2+10a_3(x+x_i)xx_i+15a_2xx_i+3a_1(x+x_i)+a_0
\end{equation}

For the generic case one must clearly reinstate the coefficients
$a_6$ and $a_5$ but this alone is not enough to render the equation
covariant which, in this case means \emph{invariant}, it being a
single relation.

The left hand side becomes invariant on dividing by $x-x_i$ since
both $$(\wp_{11},-2\wp_{12},\wp_{22})$$ and
$${\bf X^3}=\left(\frac{2xx_i}{x-x_i},-\frac{x+x_i}{x-x_i},\frac{2}{x-x_i}\right)$$ are
three dimensional representations. Note that the $x_i$ here can be
either choice from $x_1$ and $x_2.$

On the right hand side the ratio $\frac{yy_i}{(x-x_i)^3}$ is now
also seen to be invariant but $\frac{F(x,x_i)}{(x-x_i)^3}$ is not.

Note however that there is a seven dimensional representation,
\begin{eqnarray}
{\bf X^7}=\left(\frac{6}{(x-x_i)^3},-\frac{3(x+x_i)}{(x-x_i)^3},\frac{3(x^2+3xx_i+x_i^2)}{(x-x_i)^3},-\frac{(x^3+9x^2x_i+9x_i^2x+x^3)}{(x-x_i)^3}\right.,\nonumber\\
\left.\frac{3(x^2+3xx_i+x_i^2)xx_i}{(x-x_i)^3},-\frac{3(x+x_i)x^2x_i^2}{(x-x_i)^3},\frac{6x^3x_i^3}{(x-x_i)^3}\right)
\end{eqnarray}
which, when taken with the coefficients
$a_0,-a_1,a_2,-a_3,a_4,-a_5,a_6$ gives an invariant. This
modification does not alter the fundamental requirement that in the
limit $x\rightarrow x_i,y\rightarrow y_i$ the $\wp_{ij}$ are regular
but have poles of order 2 when $x\rightarrow x_i,y\rightarrow -y_i$
\cite{CF1996}. Hence our modified definition is,
\begin{equation}
\wp_{11}{\bf X^3}_2+\wp_{12}{\bf X^3}_1+\wp_{22}{\bf
X^3}_0=\frac{\tilde F(x,x_i)-yy_i}{2(x-x_i)^3}
\end{equation}
where
\begin{equation}
\frac{\tilde F(x,x_i)}{(x-x_i)^3}=a_0{\bf X^7}_6+a_1{\bf
X^7}_5+a_2{\bf X^7}_4+a_3{\bf X^7}_3+a_4{\bf X^7}_2+a_5{\bf
X^7}_1+a_6{\bf X^7}_0
\end{equation}
is a covariant ``polar'' form.

The corresponding generalizations for other genera are
straightforward and depend on constructing $2g+3$ dimensional
representations, ${\bf X}^{2g+3}$, by taking highest weight elements
$(x-x_i)^{-(g+1)}$ for $\bf e$ and applying $\bf f$ successively,
with appropriate normalizations.

Thus, for instance, for genus one we write,
\begin{equation}
\wp_{11}=\frac{\tilde F(x,x_i)-yy_i}{2(x-x_i)^2}
\end{equation}
where, using $\bf X^5$,
\begin{equation}
\tilde
F(x,x_i)=a_0+2a_1(x+x_i)+a_2(x^2+xx_i+x_i^2)+a_3(x+x_i)xx_i+a_4x^2x_i^2.
\end{equation}

The covariant polar form stands in a geometric relation to the
hyperelliptic curve $y^2-\sum_{i=0}^{2g+2}\binom{2g+2}{i}a_ix^i$ of
degree $g+2$ not shared by the traditional polar form; namely, the
curve $yy_i-\tilde F(x,x_i)=0$ of degree $g+1$ is tangent to order
$g+1$ to the hyperelliptic curve at the common point $(x_i,y_i)$.

For the calculations which follow we put the defining relations into
the convenient form
\begin{equation}
yy_i-{\bf x}^th{\bf x}_i=0
\end{equation}
where $h$ is a $(g+2)\times(g+2)$ matrix whose entries depend only
upon the $a_i$ and the $\wp_{ij}$. The $\bf x$'s are $g+2$-vectors
of monomials,e.g.
\begin{equation}
{\bf x}^t=(1,x,x^2,\ldots,x^g,x^{g+1}).
\end{equation}

\section{Differential relations in genus one}
Here we give a new, covariant treatment of the most classical case
of all: the Weierstrass $\wp$-function.

Covariance of the quartic curve
\begin{equation}
y^2=a_0+4a_1x+6a_2x^2+4a_3x^3+a_4x^4
\end{equation}
under ${\mathfrak sl}_2({\mathbb C})$ requires
\begin{eqnarray}
{\bf e}(x)&=&1\nonumber\\
{\bf e}(y)&=&0\nonumber\\
{\bf f}(x)&=&-x^2\nonumber\\
{\bf f}(y)&=&-2xy\nonumber\\
{\bf e}(a_i)&=&-(4-i)a_{i+1}\nonumber\\
{\bf f}(a_i)&=&-ia_{i-1}\nonumber
\end{eqnarray}

There is only one holomorphic differential on the curve:
$du_1=\frac{dx}{y}$ and it is clear that
\begin{eqnarray}
{\bf e}(du_1)&=&0\nonumber\\
{\bf f}(du_1)&=&0\nonumber
\end{eqnarray}
so that $\wp_{11}$, $\wp_{111}$, etc. are all invariant.

Even for this, the simplest case, it is necessary to make the Klein
definition covariant before we start by using ${\bf X}^5$ as at the
end of the last section. We apply the fundamental definition of
Klein \cite{BEL1997}, written in the form
\begin{equation}\label{K1}
yy_1-{\bf x}^th{\bf x}_1=0
\end{equation}
where ${\bf x}=(1,x,x^2)$, ${\bf x_1}=(1,x_1,x_1^2)$ but where $h$
is now the covariantly modified, three by three matrix
\begin{equation}
h=\left[
\begin{array}{ccc}
a_0&2a_1&a_2-2\wp_{11}\\
2a_1&4a_2+4\wp_{11}&2a_3\\
a_2-2\wp_{11}&2a_3&a_4
\end{array}
\right]
\end{equation}

Note that in terms of entries of $h$,
\begin{eqnarray}
y^2&=&a(x)\nonumber\\
&=&h_{33}x^4+(h_{32}+h_{23})x^3+(h_{31}+h_{22}+h_{13})x^2\nonumber\\
&&+(h_{12}+h_{21})x+h_{11}\nonumber
\end{eqnarray}
each coefficient being independent of the $\wp_{11}$ symbol.

Take the residue of (\ref{K1}) at $x=\infty$, $y={\sqrt
h_{33}}(x^2+\frac{h_{32}}{h_{33}}x+\ldots)$:
\begin{equation}\label{K10}
{\sqrt h_{33}}y_1-h_{31}-h_{32}x_1-h_{33}x_1^2=0
\end{equation}

The two index symbol $\wp_{11}$ is \cite{BEL1997} a function of $x$
and $x_1$ in the form

\begin{equation}
\wp_{11}=\wp_{11}\left(\int^xd{\bf u}+\int^{x_1}d{\bf u}\right)
\end{equation}

Hence the effect of the operator $y\partial_x=\partial_{u_1}$ etc.
on $\wp_{11}$ is
\begin{eqnarray}
y\partial_x\wp_{11}&=&\wp_{111}\\
y_1\partial_{x_1}\wp_{11}&=&\wp_{111}
\end{eqnarray}

Now apply $y\partial_x$ to the Klein relation (\ref{K1}):
\begin{equation}
yy'y_1-y{\bf x}'^th{\bf x}_1={\bf x}^t(\partial_{u_1}h){\bf x}_1
\end{equation}

Use of the defining relation allows us to replace $yy_1$ to give:
\begin{equation}
(y'{\bf x}^t-y{\bf x'}^t)h{\bf x}_1={\bf x}^t(\partial_{u_1}h){\bf
x}_1
\end{equation}

The highest order term using $y={\sqrt
h_{33}}(x^2+\frac{h_{32}}{h_{33}}x+\ldots)$, yields

\begin{eqnarray}
h_{33}(h{\bf x}_1)_2-h_{23}(h{\bf
x}_1)_3=\sqrt{h_{33}}(\partial_{u_1}h{\bf x}_1)_3
\end{eqnarray}
where we have used subscripts $(\cdot)_2$ and $(\cdot)_3$ to denote
the second and third components of a vector quantity.

Explicitly we have the identity:

\begin{eqnarray}
\left|\begin{array}{cc}
h_{12}&h_{13}\\
h_{23}&h_{33}
\end{array}\right|
+\left|\begin{array}{cc}
h_{22}&h_{23}\\
h_{32}&h_{33}
\end{array}\right|x_1
+2{\sqrt h_{33}}\wp_{111}=0&&\nonumber
\end{eqnarray}

The same identity arises if we differentiate the Klein relation with
respect to $x_1$.

So far then $y_1$ is given by a quadratic in $x_1$, linear in
$\wp_{11}$, and $\wp_{111}$ by a relation linear in $x_1$ and
$\wp_{11}$. One further relation is afforded by the fact that
$(x_1,y_1)$ lies on the curve. Using the expression (\ref{K10}) for
$y_1$ this becomes

\begin{equation}
\left|\begin{array}{cc}h_{22}&h_{23}\\h_{32}&h_{33}\end{array}\right|x_1^2+
2\left|\begin{array}{cc}h_{12}&h_{13}\\h_{32}&h_{33}\end{array}\right|x_1+
\left|\begin{array}{cc}h_{11}&h_{13}\\h_{31}&h_{33}\end{array}\right|=0.
\end{equation}

We now eliminate $x_1$ between this quadratic relation and the
preceeding linear expression for $\wp_{111}$. We obtain
\begin{equation}
\wp_{111}^2=-\frac14\left|\begin{array}{ccc}h_{11}&h_{12}&h_{13}\\h_{21}&h_{22}&h_{23}\\h_{31}&h_{32}&h_{33}\\\end{array}\right|
\end{equation}

Identifying as customary the classical $\wp$-function with
$\wp_{11}$ and its derivative, $\wp'$, with $\wp_{111}$ we have,
expanding the determinant, the equation for the $\wp$-function for
the \emph{generic} curve of genus one:
\begin{eqnarray}
\wp'^2-4\wp^3&=&-(a_0a_4-4a_1a_3+3a_2^2)\wp\\
&&-a_0a_2a_4+a_0a_3^2-2a_1a_2a_3+a_2^3+a_1^2a_4\nonumber
\end{eqnarray}

\subsection{Remarks}
\subsubsection{}
The coefficients $a_0a_4-4a_1a_3+3a_2^2$ and
$-a_0a_2a_4+a_0a_3^2-2a_1a_2a_3+a_2^3+a_1^2a_4$ are readily verified
to be invariants under the ${\mathfrak sl}_2({\mathbb C})$ action.
This is only to be expected from the classical approach. They sit
inside the two-fold and three-fold tensor products of the five
dimensional representation spanned by $\{a_0,a_1,a_2,a_3,a_4\}$.
\subsubsection{}
Specializing to the case where one branch point is moved to
$\infty$, we take $a_4=0$. By shifting $x$ we can set $a_2=0$ and by
scaling, set $a_3=1$:
\begin{equation}
\wp'^2=4\wp^3+4a_1\wp+a_0\nonumber
\end{equation}
Traditionally one associates this curve with the cubic
\[y^2=4x^3+4a_1x+a_0\]
parametrized by setting $x=\wp$ and $y=\wp'$ but we see that in fact
the origin of the factor of 4 on the left hand side is not at all
related to the value of $a_3$. It is rather an intrinsic value that
holds for the generic curve. We could of course solve the relations
obtained in the previous section to obtain $x_1$ and $y_1$ as
functions of $\wp_{11}$, $\wp_{111}$ and the $a_i$ inorder to
parametrise the generic quartic, $y^2=a(x)$. This parametrization
looks, at first sight, rather unattractive although it reduces to
the classical one when the branch point is moved to $\infty$.
\subsubsection{}
The generic differential equation for the $\wp$-function above is
actually what for higher genus would be called a quadratic identity.
Consequently the coefficients in the differential equation are
polynomial in the $a_i$ and not linear.
\subsubsection{}
Why is life more complicated for higher genus? Simply because the
$\wp_{ij}$ are now a $\frac12g(g+1)$ dimensional (not, in general,
irreducible) representation and so their relations cannot be
constructed solely from invariant quantities.

\section{Differential relations in genus two}
The fundamental definition of Klein can be modified to the form
\begin{equation}\label{K2}
yy_i-{\bf x}h{\bf x}^T_i=0
\end{equation}
for $i=1,2$, where ${\bf x}=(1,x,x^2,x^3)$, ${\bf
x_i}=(1,x_i,x_i^2,x_i^3)$ and $h$ is the covariant four by four
matrix
\begin{equation}
h=\left[
\begin{array}{cccc}
a_0&3a_1&3a_2-2\wp_{11}&a_3-2\wp_{12}\\
3a_1&9a_2+4\wp_{11}&9a_3+2\wp_{12}&3a_4-2\wp_{22}\\
3a_2-2\wp_{11}&9a_3+2\wp_{12}&9a_4+4\wp_{22}&3a_5\\
a_3-2\wp_{12}&3a_4-2\wp_{22}&3a_5&a_6
\end{array}
\right]
\end{equation}

Note that in terms of entries of $h$,
\begin{eqnarray}
y^2&=&a(x)\nonumber\\
&=&h_{44}x^6+(h_{34}+h_{43})x^5+(h_{24}+h_{33}+h_{42})x^4\nonumber\\
&&+(h_{14}+h_{23}+h_{32}+h_{41})x^3\nonumber\\
&&+(h_{13}+h_{22}+h_{31})x^2+(h_{12}+h_{21})x+h_{11}\nonumber
\end{eqnarray}
each coefficient being independent of the $\wp_{ij}$ symbols.

Take the residue of (\ref{K2}) at $x=\infty$, $y={\sqrt
h_{44}}(x^3+\frac{h_{34}}{h_{44}}x^2+\ldots)$:
\begin{equation}\label{K0}
{\sqrt h_{44}}y_1-h_{41}-h_{42}x_1-h_{43}x_1^2-h_{44}x_1^3=0
\end{equation}

The two index symbols, $\wp_{ij}$ are \cite{BEL1997} functions of
$x$, $x_1$ and $x_2$ in the form

\begin{equation}
\wp_{ij}=\wp_{ij}\left(\int^xd{\bf u}+\int^{x_1}d{\bf
u}+\int^{x_2}d{\bf u}\right)
\end{equation}

Hence the effect of the operators
$y\partial_x=\partial_{u_1}+x\partial_{u_2}$ etc. on the $\wp_{ij}$
is
\begin{eqnarray}
y\partial_x\wp_{ij}&=&\wp_{ij1}+x\wp_{ij2}\\
y_1\partial_{x_1}\wp_{ij}&=&\wp_{ij1}+x_1\wp_{ij2}\\
y_2\partial_{x_2}\wp_{ij}&=&\wp_{ij1}+x_2\wp_{ij2}
\end{eqnarray}

Apply $y_2\partial_{x_2}$ to the Klein relation (\ref{K2}) with
$i=1$. By elementary algebra it reduces, for all $x$, to the form
\begin{equation}
-2(x-x_1)^2\left(A+xB\right)=0
\end{equation}
where $A$ and $B$ are functions of $x_1$, $x_2$ and the $\wp_{ijk}$.
As there can be no relation linear in $x$ between these objects
\cite{B1897}, both the coefficients $A$ and $B$ must vanish:
\begin{eqnarray}\label{KK}
\wp_{111}+(x_1+x_2)\wp_{112}+x_1x_2\wp_{122}&=&0\\
\wp_{112}+(x_1+x_2)\wp_{122}+x_1x_2\wp_{222}&=&0\nonumber
\end{eqnarray}

Now apply $y\partial_x$ to the Klein relation (\ref{K2}) with $i=1$:
\begin{equation}
yy'y_1-y{\bf x}'h{\bf x}^T_1={\bf
x}(\partial_{u_1}h+x\partial_{u_2}){\bf x}^T_1
\end{equation}

Use of the defining relation allows us to replace $yy_1$ to give:
\begin{equation}
(y'{\bf x}^T-y{\bf x'}^T)h{\bf x}^T_1={\bf
x}(\partial_{u_1}h+x\partial_{u_2}){\bf x}^T_1
\end{equation}

Using $y={\sqrt h_{44}}(x^3+\frac{h_{34}}{h_{44}}x^2+\ldots)$, the
highest order term yields
\begin{eqnarray}
h_{44}(h{\bf x}^T_1)_3-h_{34}(h{\bf
x}^T_1)_4=\sqrt{h_{44}}(\partial_{u_2}h{\bf x}_1^T)_4
\end{eqnarray}
where again we have used subscripts $(\cdot)_i$ to denote $i$th
components of a vector quantity.

Explicitly we have a quadratic identity:

\begin{eqnarray}\label{KKK}
\left|\begin{array}{cc}
h_{31}&h_{34}\\
h_{41}&h_{44}
\end{array}\right|
+\left|\begin{array}{cc}
h_{32}&h_{34}\\
h_{42}&h_{44}
\end{array}\right|x_1
+\left|\begin{array}{cc}
h_{33}&h_{34}\\
h_{43}&h_{44}
\end{array}\right|x_1^2&&\\+2{\sqrt
h_{44}}(\wp_{122}+x_1\wp_{222})=0&&\nonumber
\end{eqnarray}

By the general symmetry of the problem the same identity must be
satisfied by $x_2$. Thus we can obtain expressions for the symmetric
combinations $x_1+x_2$ and $x_1x_2$, namely:

\begin{eqnarray}
2{\sqrt h_{44}}\wp_{222}&=&-\left|\begin{array}{cc}
h_{32}&h_{34}\\
h_{42}&h_{44}
\end{array}\right|-\left|\begin{array}{cc}
h_{33}&h_{34}\\
h_{43}&h_{44}
\end{array}\right|(x_1+x_2)\\
2{\sqrt h_{44}}\wp_{122}&=&-\left|\begin{array}{cc}
h_{31}&h_{34}\\
h_{41}&h_{44}
\end{array}\right|+\left|\begin{array}{cc}
h_{33}&h_{34}\\
h_{43}&h_{44}
\end{array}\right|x_1x_2
\end{eqnarray}

Eliminating these symmetric combinations from the second of the pair
(\ref{KK}) we obtain the relation:

\begin{eqnarray}
\left|\begin{array}{cc}
h_{33}\wp_{112}-h_{32}\wp_{122}+h_{31}\wp_{222}& h_{34}\\
&\\h_{43}\wp_{112}-h_{42}\wp_{122}+h_{41}\wp_{222}&h_{44}
\end{array}\right|&=&0\nonumber\\
&&\nonumber
\end{eqnarray}
from which it follows that
\begin{eqnarray}
h_{33}\wp_{112}-h_{32}\wp_{122}+h_{31}\wp_{222}&=&\lambda h_{34}\\
h_{43}\wp_{112}-h_{42}\wp_{122}+h_{41}\wp_{222}&=&\lambda h_{44}
\end{eqnarray}
$\lambda$ being some constant to be determined.

All the elements of these identities belong to irreducible
representations of ${\mathfrak sl}_2$ and it is easy to show that
the identities are mutually self-consistent under the Lie algebra
action if $\lambda$ is identified with $\wp_{111}$. They then form
two of a multiplet of four identities (a four dimensional
representation of ${\mathfrak sl}_2$) summarized in matrix form as,

\begin{equation}\label{bilinear}
\left(\begin{array}{cccc}
h_{11}&h_{12}&h_{13}&h_{14}\\
h_{21}&h_{22}&h_{23}&h_{24}\\
h_{31}&h_{32}&h_{33}&h_{34}\\
h_{41}&h_{42}&h_{43}&h_{44}
\end{array}\right)\left(\begin{array}{c}\wp_{222}\\-\wp_{122}\\\wp_{112}\\-\wp_{111}\end{array}\right)=0
\end{equation}

An immediate consequence of this is the relation for the Kummer
surface, quartic in the $\wp_{ij}$:
\begin{equation}
\det(h)=0
\end{equation}

But, more than this, it follows (we do not give the argument here
because it is a simplification of that leading up to equation
(\ref{seelater}) for the genus three case) straightforwardly from
(\ref{bilinear}) and the theory of diagonalisation of the symmetric
matrix $h$ that if we define the bordered matrix,

\begin{equation}
H=\left(\begin{array}{ccccc}
h_{11}&-h_{12}&h_{13}&-h_{14}&l_0\\
-h_{21}&h_{22}&-h_{23}&h_{24}&l_1\\
h_{31}&-h_{32}&h_{33}&-h_{34}&l_2\\
-h_{41}&h_{42}&-h_{43}&h_{44}&l_3\\
l_0&l_1&l_2&l_3&0
\end{array}\right)
\end{equation}
then $\det(H)$ is, up to a factor, the expression
$(l_0\wp_{222}+l_1\wp_{122}+l_2\wp_{112}+l_3\wp_{111})^2$.

That this factor is, in fact, $-\frac14$ could be established by the
classical argument of singularity balancing between the leading
terms, quadratic in the $\wp_{ijk}$ and cubic in the $\wp_{ij}$.
However it is instructive and in keeping with the current, purely
algebraic, philosophy to establish the result by using the relation
arising by application of $y_1\partial_{x_1}$ to the Klein relation
(\ref{K2}) for $i=1$.

Immediately we have
\begin{equation}
yy_1y'_1-{\bf x}(\partial_{u_1}h+x_1\partial_{u_2}h){\bf
x_1}^T-y_1{\bf x}h{\bf x'_1}^T=0
\end{equation}

Replacing $y_1y'_1$ by $\frac12a'(x_1)$, taking the $x=\infty$
residue of
\begin{equation}
\frac12ya'(x_1)-y_1{\bf x}h{\bf x'_1}^T={\bf
x}(\partial_{u_1}h+x_1\partial_{u_2}h){\bf x_1}^T
\end{equation}
and by elimination of $y_1$ as before, we find:
\begin{equation}
\frac12\left((h{\bf x}^T_1)_4^2-h_{44}a(x_1)\right)'=2{\sqrt
h_{44}}(\wp_{112}+2\wp_{122}x_1+\wp_{222}x_1^2)\nonumber
\end{equation}
prime denoting differentiation with respect to $x_1$: that is, we
ignore the implicit $x_1$ dependence of the $\wp_{ij}$. In fact the
right hand side of this equation is easily seen to be cubic in $x_1$
and not, as at first sight it appears, quintic. Exploiting the
symmetry of $h$ gives us,
\begin{eqnarray}
{\sqrt h_{44}}(\wp_{112}+2\wp_{122}x_1+\wp_{222}x_1^2)&+&
\left|\begin{array}{cc}
h_{33}&h_{34}\\
h_{43}&h_{44}
\end{array}\right|x_1^3
+\frac32\left|\begin{array}{cc}
h_{23}&h_{24}\nonumber\\
h_{43}&h_{44}
\end{array}\right|x_1^2\\
&&+ \left(\left|\begin{array}{cc}
h_{13}&h_{14}\\
h_{43}&h_{44}
\end{array}\right|+\frac12\left|\begin{array}{cc}
h_{22}&h_{24}\\
h_{42}&h_{44}
\end{array}\right|\right)x_1
\nonumber\\
&&+\frac12\left|\begin{array}{cc}
h_{12}&h_{14}\\
h_{42}&h_{44}
\end{array}\right|=0
\end{eqnarray}

It is straightforward to eliminate (\ref{KKK}) from the above to
leave a second quadratic identity,
\begin{equation}
2{\sqrt h_{44}}(\wp_{112}-\wp_{222}x_1^2)+ \left|\begin{array}{cc}
h_{23}&h_{24}\\ h_{43}&h_{44}
\end{array}\right|x_1^2
+\left|\begin{array}{cc}
h_{22}&h_{24}\\
h_{42}&h_{44}
\end{array}\right|x_1
+\left|\begin{array}{cc}
h_{12}&h_{14}\\
h_{42}&h_{44}
\end{array}\right|=0\nonumber
\end{equation}

Again, eliminating $x_1^2$ between this and (\ref{KKK}) provides a
relation of degree one in $x_1$. Since the $x_i$ can satisfy nothing
simpler than quadratic relations the coefficient of $x_1$ and the
constant term must be identically zero. The first is
\begin{equation}
4h_{44}\wp_{222}^2= \left|
\begin{array}{cc}
h_{32}&h_{34}\\
h_{42}&h_{44}
\end{array}
\right|
\left|
\begin{array}{cc}
h_{23}&h_{24}\\
h_{43}&h_{44}
\end{array}\right|
- \left|
\begin{array}{cc}
h_{33}&h_{34}\\
h_{43}&h_{44}
\end{array}\right|
\left|
\begin{array}{cc}
h_{22}&h_{24}\\
h_{42}&h_{44}
\end{array}\right|
\end{equation}
and by a well-known identity for $3\times 3$ determinants
\cite{Aitken}:
\begin{equation}
-4\wp_{222}^2= \left|
\begin{array}{ccc}
h_{22}&h_{23}&h_{24}\\
h_{32}&h_{33}&h_{34}\\
h_{42}&h_{43}&h_{44}
\end{array}\right|
\end{equation}

This allows us to fix the value of the constant of proportionality
and we obtain a beautiful, covariant generalization of Baker's
formula \cite{B1907}:
\begin{equation}
(l_0\wp_{222}+l_1\wp_{122}+l_2\wp_{112}+l_3\wp_{111})^2=-\frac14\left|\begin{array}{ccccc}
h_{11}&-h_{12}&h_{13}&-h_{14}&l_0\\
-h_{21}&h_{22}&-h_{23}&h_{24}&l_1\\
h_{31}&-h_{32}&h_{33}&-h_{34}&l_2\\
-h_{41}&h_{42}&-h_{43}&h_{44}&l_3\\
l_0&l_1&l_2&l_3&0
\end{array}\right|\nonumber
\end{equation}

For later comparison we change the sign of $l_1$ and $l_3$:
\begin{equation}\label{quadtwo}
(l_0\wp_{222}-l_1\wp_{122}+l_2\wp_{112}-l_3\wp_{111})^2=-\frac14\left|\begin{array}{ccccc}
h_{11}&h_{12}&h_{13}&h_{14}&l_0\\
h_{21}&h_{22}&h_{23}&h_{24}&l_1\\
h_{31}&h_{32}&h_{33}&h_{34}&l_2\\
h_{41}&h_{42}&h_{43}&h_{44}&l_3\\
l_0&l_1&l_2&l_3&0
\end{array}\right|\nonumber
\end{equation}

In section \ref{fourindex} we will derive identities linear in the
$\wp_{ij}$ and $\wp_{ijk}$ from the above quadratic identities.
\emph{Presumably} all $\wp$-function identities arise from these
quadratic ones by algebraic and differential processes but, of
course, this is not immediately clear. Nor is it immediately
essential to their application in this paper.

\section{Differential relations in genus three.}

The last section recovers classical results in that the covariant
identities were written down in \cite{B1907} though not there
derived in a covariant manner. By contrast a covariant treatment of
higher genus hyperelliptic (or non-hyperelliptic) curves has not
been given before. This we now do.

For genus three we have three covariant Klein equations:
\begin{equation}\label{K3}
yy_i-{\bf x}h{\bf x}^T_i=0
\end{equation}
for $i=1,2,3$, where ${\bf x}=(1,x,x^2,x^3,x^4)$, ${\bf
x_i}=(1,x_i,x_i^2,x_i^3,x_i^4)$ and $h$ is the $5\times 5$ matrix,
\begin{tiny}
\begin{equation}
\left[\begin{array}{ccccc}
a_0&4a_1&6a_2-2\wp_{11}&4a_3-2\wp_{12}&a_4-2\wp_{13}\\
4a_1&16a_2+4\wp_{11}&24a_3+2\wp_{12}&16a_4-2\wp_{22}+4\wp_{13}&4a_5-2\wp_{23}\\
6a_2-2\wp_{11}&24a_3+2\wp_{12}&36a_4+4\wp_{22}-4\wp_{13}&24a_5+2\wp_{23}&6a_6-2\wp_{33}\\
4a_3-2\wp_{12}&16a_4-2\wp_{22}+4\wp_{13}&24a_5+2\wp_{23}&16a_6+4\wp_{33}&4a_7\\
a_4-2\wp_{13}&4a_5-2\wp_{23}&6a_6-2\wp_{33}&4a_7&a_8
\end{array}
\right]\nonumber
\end{equation}
\end{tiny}

The residue of (\ref{K3}) at $x=\infty$, $y(x)={\sqrt
h_{55}}x^4+\frac{h_{45}+h_{54}}{2{\sqrt h_{55}}}x^3\ldots$ gives
\begin{equation}
{\sqrt h_{5,5}}y_i-(h{\bf x}^T_i)_5=0,
\end{equation}
for $i$ with value 1, 2 or 3.

This time the operator $y\partial_x$ and its indexed relatives is
given by $\partial_{u_1}+x\partial_{u_2}+x^2\partial_{u_3}$ etc. We
may apply $y_2\partial_{x_2}$ to (\ref{K3}) with $i=1$ and take the
residue at $x=\infty$ to obtain,
\begin{eqnarray}
\left((\frac{\partial h}{\partial u_1}+x_2\frac{\partial h}{\partial
u_2}+x_2^2\frac{\partial h}{\partial u_3}){\bf x}^T_1\right)_5=0
\end{eqnarray}

Simplifying, removing overall factors of $x_1-x_2$, gives
\begin{eqnarray}
\wp_{113}+(x_1+x_2)\wp_{123}+x_1x_2\wp_{223}+(x_1^2+x_2^2)\wp_{133}&&\nonumber\\
+(x_1+x_2)x_1x_2\wp_{233}+x_1^2x_2^2\wp_{333}&=&0
\end{eqnarray}
and, by cyclic interchange of the $x_i$,
\begin{eqnarray}
\wp_{113}+(x_2+x_3)\wp_{123}+x_2x_3\wp_{223}+(x_2^2+x_3^2)\wp_{133}&&\nonumber\\
+(x_2+x_3)x_2x_3\wp_{233}+x_2^2x_3^2\wp_{333}&=&0\\
\wp_{113}+(x_3+x_1)\wp_{123}+x_3x_1\wp_{223}+(x_3^2+x_1^2)\wp_{133}&&\nonumber\\
+(x_3+x_1)x_3x_1\wp_{233}+x_3^2x_1^2\wp_{333}&=&0.
\end{eqnarray}

From these three identities we can form three identities whose
coefficients are symmetric functions in the $x_i$, namely:
\begin{eqnarray}\label{A}
\wp_{223}-\wp_{133}+s^{(1)}\wp_{233}+s^{(2)}\wp_{333}&=&0\\
\wp_{123}+s^{(1)}\wp_{133}-s^{(3)}\wp_{333}=0\\
\wp_{113}-s^{(2)}\wp_{133}-s^{(3)}\wp_{233}=0
\end{eqnarray}
where $s^{(1)}=x_1+x_2+x_3$, $s^{(2)}=x_1x_2+x_2x_3+x_3x_1$ and
$s^{(3)}=x_1x_2x_3$.

An important observation at this point is that these three equations
are overdetermined for $s^{(1)}$, $s^{(2)}$ and $s^{(3)}$ so that
the $\wp_{ijk}$ must satisfy the quadratic identity
\begin{equation}
\wp_{113}\wp_{333}-\wp_{123}\wp_{233}+\wp_{223}\wp_{133}-\wp_{133}^2=0.
\end{equation}
This relation is in the kernel of $\bf e$ and thus is a highest
weight element for a set of relations forming a five dimensional
representation:
\begin{eqnarray}\label{trivial}
P_{\bf
5}(0)&=&\wp_{113}\wp_{333}-\wp_{123}\wp_{233}+\wp_{223}\wp_{133}-\wp_{133}^2\nonumber\\
P_{\bf
5}(1)&=&-\wp_{233}\wp_{113}-\wp_{112}\wp_{333}-\wp_{133}\wp_{222}+2\wp_{133}\wp_{123}+\wp_{233}\wp_{122}\nonumber\\
P_{\bf
5}(2)&=&\wp_{133}\wp_{122}-\wp_{133}\wp_{113}-\wp_{223}\wp_{122}+\wp_{223}\wp_{113}\nonumber\\
&&+\wp_{111}\wp_{333}+\wp_{123}\wp_{222}-2\wp_{123}^2\nonumber\\
P_{\bf
5}(3)&=&-\wp_{233}\wp_{111}-\wp_{112}\wp_{133}+\wp_{112}\wp_{223}-\wp_{113}\wp_{222}+2\wp_{113}\wp_{123}\nonumber\\
P_{\bf
5}(4)&=&-\wp_{123}\wp_{112}+\wp_{113}\wp_{122}-\wp_{113}^2+\wp_{133}\wp_{111}\nonumber
\end{eqnarray}
This gives a set of five identities quadratic in the $\wp_{ijk}$,
$P_{\bf 5}(i)=0$ for $i=0,\dots 4.$

Differentiating (\ref{K3}) with respect to $y\partial_x$,
\begin{equation}
(y'(x){\bf x}-y(x){\bf x}')h{\bf x}_1^T={\bf x}(\frac{\partial
h}{\partial u_1}+x\frac{\partial h}{\partial u_2}+x^2\frac{\partial
h}{\partial u_3}){\bf x}_1^T
\end{equation}
we again use the expansion near $x=\infty$,
\[y(x)={\sqrt h_{55}}x^4+\frac{h_{45}+h_{54}}{2{\sqrt
h_{55}}}x^3\ldots\] collecting the highest order term (degree 6) in
the identity just obtained. Thus, using the symmetry of the matrix
$h$,
\begin{equation}\label{cubic}
\left|\begin{array}{cc}
(h{\bf x}_1^T)_4 & (h{\bf x}_1^T)_5\\
h_{54} & h_{55}
\end{array}\right|=(\frac{\partial
h}{\partial u_3}{\bf x}_1^T)_5.
\end{equation}

This is an identity cubic in $x_1$ also satisfied by $x_2$ and
$x_3$:
\begin{eqnarray}
\left|\begin{array}{cc}
h_{44} & h_{45}\\
h_{54} & h_{55}
\end{array}\right|x_i^3+
\left|\begin{array}{cc}
h_{34} & h_{35}\\
h_{54} & h_{55}
\end{array}\right|x_i^2+
\left|\begin{array}{cc}
h_{24} & h_{25}\\
h_{44} & h_{55}
\end{array}\right|x_i+
\left|\begin{array}{cc}
h_{14} & h_{15}\\
h_{54} & h_{55}
\end{array}\right|&&\nonumber\\
+2{\sqrt h_{55}}(\wp_{133}+x_i\wp_{233}+x_i^2\wp_{333})&=&0
\end{eqnarray}

We can now eliminate the symmetric functions $s^{(1)}$, $s^{(2)}$
and $s^{(3)}$ in (\ref{A}). From the first relation we obtain,
\begin{eqnarray}
h_{24}\wp_{333}-h_{34}\wp_{233}+h_{44}(\wp_{223}-\wp_{133})-
h_{54}\lambda&=&0\nonumber\\
h_{25}\wp_{333}-h_{35}\wp_{233}+h_{45}(\wp_{223}-\wp_{133})-
h_{55}\lambda&=&0
\end{eqnarray}
where $\lambda$ is an undetermined multiplier.

Since these identities are polynomial in the $\wp_{ijk}$ and the
$h_{ij}$ only they must belong to a finite dimensional
representation of ${\mathfrak sl}_2(\mathbb C)$. Application of the
$\bf e$ and $\bf f$ operators must generate further identities. This
can only work for a special value of $\lambda$ and application of
$\bf e$ to the second of the above identities shows that it will be
a highest weight element (in the kernel of $\bf e$) only if
$\lambda=\wp_{222}-2\wp_{123}.$ Hence we have
\begin{eqnarray}\label{base1}
h_{24}\wp_{333}-h_{34}\wp_{233}+h_{44}(\wp_{223}-\wp_{133})-h_{54}(\wp_{222}-2\wp_{123})
&=&0\nonumber\\
h_{25}\wp_{333}-h_{35}\wp_{233}+h_{45}(\wp_{223}-\wp_{133})-h_{55}(\wp_{222}-2\wp_{123})
&=&0
\end{eqnarray}

We label the second of these identities $P_{\bf 9}(0)$ because it is
highest weight for a nine dimensional representation generated by
repeated application of $\bf f$, a set of nine linearly independent
identities $P_{\bf 9}(i)$ for $i=0,\ldots 8$. The last of these
identities is:
\begin{equation}
P_{\bf
9}(8)=h_{11}(\wp_{222}-2\wp_{123})-h_{12}(\wp_{122}-\wp_{113})+h_{13}\wp_{112}-h_{14}\wp_{111}=0
\end{equation}

Rather than write these out in detail now we shall summarize them in
a more compact form shortly.

Now from a linear combination of the first of the identities
(\ref{base1}) and $P_{\bf 9}(1)$ we can form the highest weight
identity for a seven dimensional representation:
\begin{eqnarray}
P_{\bf
7}(0)=-4h_{15}\wp_{333}+4h_{35}\wp_{133}-h_{45}(2\wp_{123}+\wp_{222})+4h_{55}(\wp_{122}-\wp_{113})&&\nonumber\\
-h_{34}\wp_{233}+h_{24}\wp_{333}-h_{44}(\wp_{133}-\wp_{223})&=&0
\end{eqnarray}

Proceeding in this way with the other identities obtained from
eliminating the symmetric functions from the other identities in
(\ref{A}), we obtain representations $P_{\bf 5}$, $P_{\bf 3}$ and
$P_{\bf 1}$, giving a total of $9+7+5+3+1=5^2$ relations linear in
the three index symbols.

These identities are not presented in the simplest form however.
They can be rendered more transparent by taking various linear
combinations so that one only ever has four $h$ terms arising in
each identity. We do not give the details here because it involves
routine linear algebra applied to the above identities, best
accomplished using a computer algebra package. The \emph{fact} that
this simplification is possible, however, is of significance and it
not clear to the present author exactly why it should be so.

After this rearrangement the identities take the form of a matrix
product
\begin{equation}\label{fourterm}
hA=0 \end{equation}
of the symmetric $5\times 5$ matrix $h$ and an
antisymmetric $5\times 5$ matrix

\begin{equation}
A=\left[
\begin{array}{ccccc}
0&-\wp_{333}&\wp_{233}&-\wp_{223}+\wp_{133}&\wp_{222}-2\wp_{123}\\
\wp_{333}&0&-\wp_{133}&\wp_{123}&-\wp_{122}+\wp_{113}\\
-\wp_{233}&\wp_{133}&0&-\wp_{113}&\wp_{112}\\
\wp_{223}-\wp_{133}&-\wp_{123}&\wp_{113}&0&-\wp_{111}\\
-\wp_{222}+2\wp_{123}&\wp_{122}-\wp_{113}&-\wp_{112}&\wp_{111}&0
\end{array}
\right]
\end{equation}

One checks that the matrix $A$ has rank at most three by virtue of
the relations (\ref{trivial}) obtained earlier. In fact the $4\times
4$ minors of $A$ are products of the $P_{\bf 5}(i)$:
\begin{equation}
A(i,j)=P_{\bf 5}(5-i)P_{\bf 5}(5-j)
\end{equation}
Further the $3\times 3$ minors also have the $P_{\bf 5}(i)$ as
factors. There are however non vanishing $2\times 2$ minors so the
rank of $A$ is exactly two.

Consequently the five by five matrix $h$ has exactly a two
dimensional zero eigenspace and, being symmetric, must be similar to
a diagonal matrix of form $h_D=Diag(0,0,h_3,h_4,h_5)$.

We can use this fact to obtain identities quadratic in the
$\wp_{ijk}$ by generalizing the argument in the genus two case as
follows.

Let $\Pi$ be the matrix which diagonalizes $h$, let $\bf l$ and $\bf
k$ be arbitrary five component column vectors, $I_2$ the two by two
identity matrix and consider
\begin{eqnarray}
\left|\left[\begin{array}{cc} \Pi^T & 0\\0 & I_2
\end{array}\right]\left[\begin{array}{ccc}
h&\bf l & \bf
k\\
\bf l^T & 0 & 0\\
\bf k^T & 0 & 0\\\end{array}\right]\left[\begin{array}{cc}\Pi & 0\\0
&
I_2\end{array}\right]\right|&=&\left|\left[\begin{array}{ccc}h_D&\Pi^T\bf
l & \Pi^T\bf
k\\
{\bf l^T}\Pi & 0 & 0\\
{\bf k^T}\Pi & 0 & 0\\\end{array}\right]\right|\nonumber\\
&=&h_3h_4h_5\left|\left[\begin{array}{cc}(\Pi^T{\bf l})_1 &
(\Pi^T{\bf l})_2\\(\Pi^T{\bf k})_1 & (\Pi^T{\bf k})_2\end{array}
\right]\right|^2
\end{eqnarray}

Now consider
\begin{eqnarray}{\bf l}^TA{\bf k}&=&{\bf l}^T\Pi A_D \Pi^T{\bf
k}\nonumber\\
&=&\alpha\left|\left[\begin{array}{cc}(\Pi^T{\bf l})_1 & (\Pi^T{\bf
l})_2\\(\Pi^T{\bf k})_1 & (\Pi^T{\bf k})_2\end{array} \right]\right|
\end{eqnarray}
where $A_D$ is the normal form of $A$
\begin{equation}
A_D=\left[\begin{array}{ccccc}
0 & \alpha & 0 & 0 & 0\\
-\alpha & 0 & 0 & 0 & 0\\
0 & 0 & 0 & 0 & 0\\
0 & 0 & 0 & 0 & 0\\
0 & 0 & 0 & 0 & 0\\
\end{array}\right]
\end{equation}
corresponding to the diagonal form of $h$.

Combining these observations we obtain the attractive formula
\begin{equation}\label{seelater}
({\bf l}^TA{\bf k})^2=\lambda\left|\begin{array}{ccc} h&\bf l & \bf
k\\
\bf l^T & 0 & 0\\
\bf k^T & 0 & 0\\\end{array}\right|
\end{equation}
where $\lambda$ is a function yet to be determined.

The undetermined factor can be found from a (simple) singularity
argument and also by a more involved, algebraic expansion and as for
genus two we will present the latter.

We return to the original relation on the curve for $y_1$ of degree
8 in $x_1$. By substituting for $y_1$ we obtain a sextic in $x_1$
from which we eliminate the degree 6 and 5 terms using the cubic
expression (\ref{cubic}). The resulting quartic identity in $x_1$
has leading term
\begin{equation}
\wp_{333}^2+\frac14\left|\begin{array}{ccc}
h_{33}&h_{34}&h_{35}\\
h_{43}&h_{44}&h_{45}\\
h_{53}&h_{54}&h_{55}
\end{array}\right|.
\end{equation}

Hence $\lambda=-\frac14$ in the full quadratic identity:
\begin{equation}
({\bf l}^TA{\bf k})^2=-\frac14\left|\begin{array}{ccc} h&\bf l & \bf
k\\
\bf l^T & 0 & 0\\
\bf k^T & 0 & 0\\\end{array}\right|
\end{equation}

This formula is a new result of this paper.

\section{Identities for $\wp_{ijkl}$}\label{fourindex}
In all three cases discussed above there are identities for the four
index $\wp$-functions of the form:
\begin{equation}
\wp_{ijkl}=F(\wp_{11},\wp_{12},\wp_{22},\wp_{13},\ldots)
\end{equation}
that are obtained by differentiating the identities quadratic in the
$\wp_{ijk}$ and (for genus two and three) using certain identities
involving two and three index $\wp$-functions.

Clearly in genus one we get
\begin{equation}
\wp''=6\wp^2-\frac12(a_0a_4-4a_1a_3+3a_2^2)
\end{equation}
recalling that the two index $\wp$ function is written as $\wp$ in
this case.

In genus two we start with the identity for $\wp_{222}^2$.
Differentiating,
\begin{eqnarray}
-8\wp_{222}\wp_{2222}&=& \left|
\begin{array}{ccc}
4\wp_{112}&h_{23}&h_{24}\\
2\wp_{122}&h_{33}&h_{34}\\
-2\wp_{222}&h_{43}&h_{44}
\end{array}\right|+
\left|
\begin{array}{ccc}
h_{22}&2\wp_{122}&h_{24}\\
h_{32}&4\wp_{222}&h_{34}\\
h_{42}&0&h_{44}
\end{array}\right|\nonumber\\
&&+ \left|
\begin{array}{ccc}
h_{22}&h_{23}&-2\wp_{222}\\
h_{32}&h_{33}&0\\
h_{42}&h_{43}&0
\end{array}\right|\nonumber\\
&=&4\wp_{112}\left|\begin{array}{cc}h_{33}&h_{34}\\h_{43}&h_{44}\end{array}\right|
+4\wp_{122}\left|\begin{array}{cc}h_{32}&h_{34}\\h_{42}&h_{44}\end{array}\right|\nonumber\\
&&+4\wp_{222}\left(-\left|\begin{array}{cc}h_{23}&h_{24}\\h_{33}&h_{34}\end{array}\right|+\left|\begin{array}{cc}h_{22}&h_{24}\\h_{42}&h_{44}\end{array}\right|\right)\nonumber
\end{eqnarray}

The first two terms on the right hand side can be replaced by a
single term with factor $\wp_{222}$ by utilizing the identity
\begin{equation}
\wp_{112}\left|\begin{array}{cc}h_{3,3}&h_{3,4}\\h_{4,3}&h_{4,4}\end{array}\right|
+\wp_{122}\left|\begin{array}{cc}h_{3,2}&h_{3,4}\\h_{4,2}&h_{4,4}\end{array}\right|
+\wp_{222}\left|\begin{array}{cc}h_{3,1}&h_{3,4}\\h_{4,1}&h_{4,4}\end{array}\right|=0\nonumber
\end{equation}

This, in turn, is obtained from the quadratic identity by setting
$l_i=h_{i+1,j}$ to get four identities of the form

\begin{equation}
h_{1,j}\wp_{222}-h_{2,j}\wp_{122}+h_{3,j}\wp_{112}-h_{4,j}\wp_{111}=0
\end{equation}
and eliminating $\wp_{111}$ from the pair $j=3,4$. Thus
\begin{eqnarray}
-\wp_{2222}&=&\frac12\left(-\left|\begin{array}{cc}h_{2,3}&h_{2,4}\\h_{3,3}&h_{3,4}\end{array}\right|+\left|\begin{array}{cc}h_{2,2}&h_{2,4}\\h_{4,2}&h_{4,4}\end{array}\right|-\left|\begin{array}{cc}h_{3,1}&h_{3,4}\\h_{4,1}&h_{4,4}\end{array}\right|\right)\\
\frac13(-\wp_{2222}+6\wp^2_{22})&=&a_2a_6-4a_3a_5+3a_4^2+a_6\wp_{11}-2a_5\wp_{12}+a_4\wp_{22}\nonumber
\end{eqnarray}
Application of $\bf e$ and $\bf f$ to this identity shows that it is
highest weight for a five dimensional representation reproducing the
classic partial differential equations of Baker \cite{B1907}.

The fully general, covariant genus three equations have not been
written down before. We proceed as before by differentiating the
$\wp_{333}^2$ relation with respect to $u_3$:
\begin{eqnarray}
-8\wp_{333}\wp_{3333}&=&\left|\begin{array}{ccc}
4\wp_{223}-4\wp_{133}&h_{34}&h_{35}\\
2\wp_{233}&h_{44}&h_{45}\\
-2\wp_{333}&h_{54}&h_{55}
\end{array}\right|+\left|\begin{array}{ccc}
h_{33}&2\wp_{233}&h_{35}\\
h_{43}&4\wp_{333}&h_{45}\\
h_{53}&0&h_{55}
\end{array}\right|\nonumber\\
&&+\left|\begin{array}{ccc}
h_{33}&h_{34}&-2\wp_{333}\\
h_{43}&h_{44}&0\\
h_{53}&h_{54}&0
\end{array}\right|\nonumber\\
&=&4(\wp_{223}-\wp_{133})\left|\begin{array}{cc}h_{44}&h_{45}\\h_{54}&h_{55}\end{array}\right|
-4\wp_{233}\left|\begin{array}{cc}h_{34}&h_{35}\\h_{54}&h_{55}\end{array}\right|\nonumber\\
&&+4\wp_{333}\left(\left|\begin{array}{cc}h_{33}&h_{35}\\h_{53}&h_{55}\end{array}\right|-\left|\begin{array}{cc}h_{34}&h_{35}\\h_{44}&h_{45}\end{array}\right|\right)\nonumber
\end{eqnarray}

As before we can derive identities linear in the $\wp_{ijk}$ from
the general quadratic identities. Putting
$k_0=1,\,k_1=k_2=k_3=k_4=0$ and $l_i=h_{i+1,j},\,i=0\ldots 4$ gives,
for any choice of $j$,
\begin{equation}
h_{2j}\wp_{333}-h_{3j}\wp_{233}+h_{4j}(\wp_{223}-\wp_{133})-h_{5j}(\wp_{222}-2\wp_{123})=0
\end{equation}

Eliminating the $\wp_{222}-2\wp_{123}$ term between the cases $j=4$
and $j=5$ yields
\begin{equation}
(\wp_{223}-\wp_{133})\left|\begin{array}{cc}h_{44}&h_{45}\\h_{54}&h_{55}\end{array}\right|
-\wp_{233}\left|\begin{array}{cc}h_{34}&h_{35}\\h_{54}&h_{55}\end{array}\right|
+\wp_{333}\left|\begin{array}{cc}h_{24}&h_{25}\\h_{54}&h_{55}\end{array}\right|=0
\end{equation}

Use of this identity in the equation for $\wp_{333}\wp_{3333}$ gives
\begin{equation}
-2\wp_{3333}=-\left|\begin{array}{cc}h_{24}&h_{25}\\h_{54}&h_{55}\end{array}\right|
+\left|\begin{array}{cc}h_{33}&h_{35}\\h_{53}&h_{55}\end{array}\right|
-\left|\begin{array}{cc}h_{34}&h_{35}\\h_{44}&h_{45}\end{array}\right|
\end{equation}
and thus, by application of $\bf f$, the nine dimensional space of
identities,
\begin{eqnarray}
-\wp_{3333}+6\wp^2_{33}&=&10(a_4a_8-4a_5a_7+3a_6^2)\nonumber\\
&&+8a_6\wp_{33}-8a_7wp_{23}+a_8(3\wp_{22}-4\wp_{13})\nonumber\\
-\wp_{2333}+6\wp_{23}\wp_{33}&=&10(a_3a_8-3a_4a_7+2a_5a_6)\nonumber\\
&&+12a_5\wp_{33}-10\wp_{23}+4a_7(\wp_{22}-3\wp_{13})\nonumber\\
&&+2a_8\wp_{12}\nonumber\\
2(-\wp_{1333}+6\wp_{13}\wp_{33})&=&10(3a_2a_8-4a_3a_7-11a_4a_6+12a_5^2)\nonumber\\
+3(-\wp_{2233}+2\wp_{22}\wp_{33}+4\wp_{23}^2)&&+60a_4\wp_{33}-36a_5\wp_{23}-2a_6(9\wp_{22}-52\wp_{13})\nonumber\\
&&+20a_7\wp_{12}+4a_8\wp_{11}\nonumber\\
-\wp_{2223}+6\wp_{22}\wp_{23}&=&10(a_1a_8+2a_2a_7-12a_3a_6+9a_4a_5)\nonumber\\
+3(-\wp_{1233}+2\wp_{12}\wp_{33}+4\wp_{13}\wp_{23})&&+40a_3\wp_{33}-10a_4\wp_{23}+4a_5(3\wp_{22}-29\wp_{13})\nonumber\\
&&+18a_6\wp_{12}+12a_7\wp_{11}\nonumber\\
-\wp_{2222}+6\wp_{22}^2&=&10(a_0a_8+12a_1a_7-22a_2a_6-36a_3a_5+45a_4^2\nonumber\\
+6(-\wp_{1133}+2\wp_{11}\wp_{33}+4\wp_{13}^2)&&+120a_2\wp_{33}+40a_3\wp_{23}+50a_4(\wp_{22}-12\wp_{13})\nonumber\\
+12(-\wp_{1223}+4\wp_{12}\wp_{23}+2\wp_{13}\wp_{22})&&+40a_5\wp_{12}+120a_6\wp_{11}\nonumber\\
-\wp_{1222}+6\wp_{12}\wp_{22}&=&10(a_0a_7+2a_1a_6-12a_2a_5+9a_3a_4)\nonumber\\
+3(-\wp_{1123}+4\wp_{12}\wp_{13}+2\wp_{11}\wp_{23})&&+12a_1\wp_{33}+18a_2\wp_{23}+4a_3(3\wp_{22}-29\wp_{13})\nonumber\\
&&-10a_4\wp_{12}+40a_5\wp_{11}\nonumber\\
2(-\wp_{1113}+6\wp_{11}\wp_{13})&=&10(3a_0a_6-4a_1a_5-11a_2a_4+12a_3^2)\nonumber\\
+3(-\wp_{1122}+2\wp_{11}\wp_{22}+4\wp_{12}^2)&&+4a_0\wp_{33}+20a_1\wp_{23}+2a_2(9\wp_{22}-52\wp_{13})\nonumber\\
&&-36a_3\wp_{12}+60a_4\wp_{11}\nonumber\\
-\wp_{1112}+6\wp_{11}\wp_{12}&=&10(a_0a_5-3a_1a_4+2a_2a_3)\nonumber\\
&&+2a_0\wp_{23}+4a_1(\wp_{22}-3\wp_{13})-10a_2\wp_{12}\nonumber\\
&&+12a_3\wp_{11}\nonumber\\
-\wp_{1111}+6\wp_{11}^2&=&10(a_0a_4-4a_1a_3+3a_2^2)\nonumber\\
&&+a_0(3\wp_{22}-4\wp_{13})-8a_1\wp_{12}+8a_2\wp_{11}\nonumber
\end{eqnarray}

Given that there are fifteen of the symbols $\wp_{ijkl}$ we expect
to be able to find a further six identities.

Thus, returning to the $\wp_{333}^2$ identity and differentiating
with respect to $u_1$ this time yields
\begin{eqnarray}
-8\wp_{333}\wp_{1333}&=&\left|\begin{array}{ccc}
4\wp_{122}-4\wp_{113}&h_{34}&h_{35}\\
2\wp_{123}&h_{44}&h_{45}\\
-2\wp_{133}&h_{54}&h_{55}
\end{array}\right|+\left|\begin{array}{ccc}
h_{33}&2\wp_{123}&h_{35}\\
h_{43}&4\wp_{133}&h_{45}\\
h_{53}&0&h_{55}
\end{array}\right|\nonumber\\
&&+\left|\begin{array}{ccc}
h_{33}&h_{34}&-2\wp_{333}\\
h_{43}&h_{44}&0\\
h_{53}&h_{54}&0
\end{array}\right|\nonumber\\
&=&4(\wp_{122}-\wp_{113})\left|\begin{array}{cc}h_{44}&h_{45}\\h_{54}&h_{55}\end{array}\right|
-4\wp_{123}\left|\begin{array}{cc}h_{34}&h_{35}\\h_{54}&h_{55}\end{array}\right|\nonumber\\
&&+4\wp_{133}\left(\left|\begin{array}{cc}h_{33}&h_{35}\\h_{53}&h_{55}\end{array}\right|-\left|\begin{array}{cc}h_{34}&h_{35}\\h_{44}&h_{45}\end{array}\right|\right)\nonumber
\end{eqnarray}

This time some appropriate identities arise by choosing
$k_0=0,k_1=1,k_2=0,k_3=0$ and $k_4=0$ and the $l_i$ as before:
\begin{equation}
h_{ij}\wp_{333}-h_{3j}\wp_{133}+h_{4j}\wp_{123}-h_{5j}(\wp_{122}-\wp_{113})=0
\end{equation}
These allow us to replace the terms on the right hand side of the
$\wp_{1333}$ equation by terms involving $\wp_{333}$ and so factor
this out to leave
\begin{equation}
-2\wp_{1333}=-\left|\begin{array}{cc}h_{14}&h_{15}\\h_{44}&h_{45}\end{array}\right|+\left|\begin{array}{cc}h_{13}&h_{53}\\h_{15}&h_{55}\end{array}\right|
\end{equation}
Applying $\bf e$ to this identity gives the $\wp_{2333}$ identity
found above. Successive applications of $\bf f$ however yield a set
of seven identities:
\begin{eqnarray}
-\wp_{1333}+6\wp_{13}\wp_{33}&=&3a_2a_8-8a_3a_7+5a_4a_6\nonumber\\
&&+3a_4\wp_{33}-10a_6\wp_{13}+4a_7\wp_{12}-a_8\wp_{11}\nonumber\\
-\wp_{1233}+2\wp_{12}\wp_{33}+4\wp_{13}\wp_{23}&=&2a_1a_8-12a_3a_6+10a_4a_5\nonumber\\
&&+4a_3\wp_{33}+2a_4\wp_{23}-20a_5\wp_{13}+6a_6\wp_{12}\nonumber\\
-\wp_{1133}+2\wp_{11}\wp_{33}+4\wp_{13}^2&=&a_0a_8+8a_1a_7-18a_2a_6-16a_3a_5+25a_4^2\nonumber\\
-\wp_{1223}+2\wp_{13}\wp_{22}+4\wp_{12}\wp_{23}&&+6a_2\wp_{33}+8a_3\wp{23}+a_4(\wp_{22}-48\wp_{13})\nonumber\\
&&+8a_5\wp_{12}+6a_6\wp_{11}\nonumber\\
-\wp_{1222}+6\wp_{12}\wp_{22}&=&16a_0a_7+20a_1a_6-156a_2a_5+120a_3a_4\nonumber\\
+6(-\wp_{1123}+2\wp_{11}\wp_{23}+4\wp_{12}\wp_{13})&&+12a_1\wp_{33}+36a_2\wp_{23}+4(3a_3\wp_{22}-44\wp_{13})\nonumber\\
&&-4a_4\wp_{12}+52a_5\wp_{11}\nonumber\\
-\wp_{1113}+6\wp_{11}\wp_{12}&=&11a_0a_6-16a_1a_5-35a_2a_4+40a_3^2\nonumber\\
-\wp_{1122}+2\wp_{11}\wp_{22}+4\wp_{12}^2&&+a_0\wp_{33}+8a_1\wp_{23}+2a_2(3\wp_{22}-19\wp_{13})\nonumber\\
&&-12a_3\wp_{12}+21a_4\wp_{11}\nonumber\\
-\wp_{1112}+6\wp_{11}\wp_{12}&=&10(a_0a_5-3a_1a_4+2a_2a_3)\nonumber\\
&&+2a_0\wp_{23}+4a_1(\wp_{22}-3\wp_{13})-10a_2\wp_{12}\nonumber\\
&&+12a_3\wp_{11}\nonumber\\
-\wp_{1111}+6\wp_{11}^2&=&10(a_0a_4-4a_1a_3+3a_2^2)\nonumber\\
&&+a_0(3\wp_{22}-4\wp_{13})-8a_1\wp_{12}+8a_2\wp_{11}\nonumber
\end{eqnarray}

Of these seven the last two are already represented in the previous
set so that we still seek another one. To find this go to the
quadratic identity for $\wp_{133}$,
\begin{equation}
-\wp_{133}^2=\frac14\left|\begin{array}{ccc}h_{11}&h_{14}&h_{15}\\h_{41}&h_{44}&h_{55}\\h_{51}&h_{54}&h_{55}\end{array}\right|
\end{equation}
and differentiate with respect to $u_1$. Using similar identities to
before we find
\begin{equation}
-2\wp_{1133}=\left|\begin{array}{cc}h_{1,1}&h_{1,5}\\h_{5,1}&h_{5,5}\end{array}\right|-
\left|\begin{array}{cc}h_{1,4}&h_{1,5}\\h_{2,4}&h_{2,5}\end{array}\right|
\end{equation}
and applying $\bf f$ successively arrive at:
\begin{eqnarray}
2(-\wp_{1133}+6\wp_{13}^2)+4(\wp_{23}\wp_{12}-\wp_{13}\wp_{22})&=&a_0a_8-16a_3a_5+15a_4^2\nonumber\\
&&+8a_3\wp_{23}-2a_4(\wp_{22}+12\wp_{13})+8a_5\wp_{12}\nonumber\\
-\wp_{1123}+4\wp_{12}\wp_{13}+2\wp_{23}\wp_{11}&=&2a_0a_7-12a_2a_5+10a_3a_4\nonumber\\
&&+6a_2\wp_{23}-20a_3\wp_{13}+2a_4\wp_{12}+4a_5\wp_{11}\nonumber\\
-\wp_{1122}+2\wp_{11}\wp_{22}+4\wp_{12}^2&=&14a_0a_6-24a_1a_5-30a_2a_4+40a_3^3\nonumber\\
+2(-\wp_{1113}+6\wp_{11}\wp_{13})&&+12a_1\wp_{23}+6a_2(\wp_{22}-8\wp_{13})\nonumber\\
&&-12a_3\wp_{12}+24a_4\wp_{11}\nonumber\\
-\wp_{1112}+6\wp_{11}\wp_{12}&=&10a_0a_5-30a_1a_4+20a_2a_3\nonumber\\
&&+2a_0\wp_{23}+4a_1(\wp_{22}-3\wp_{13})\nonumber\\
&&-10a_2\wp_{12}+12a_3\wp_{11}\nonumber\\
-\wp_{1111}+6\wp_{11}^2&=&10a_0a_4-40a_1a_3+30a_2^2\nonumber\\
&&+a_0(3\wp_{22}-4\wp_{13})-8a_1\wp_{12}+8a_2\wp_{11}\nonumber
\end{eqnarray}

Only one of these is linearly independent of the identities we
already have.

In Appendix 1 we summarize these identities and in Appendix 2 we
compare them with the original, non-covariant identities of Baker
\cite{B1903}, showing that they are equivalent under a simple
transformation. To this end the identities in Appendix 1 are written
in a Baker friendly form where each involves but one of the four
index objects. This is not ideal from the representation theoretic
viewpoint however, as the identities then do not fall naturally into
multiplets.

\section{Conclusions}
This paper establishes that the use of covariant methods for
hyperelliptic curves is a practical tool in the construction and
understanding of the partial differential equations satisfied by the
$\wp$-function. In order to do so the definition of the $\wp$
function has to be slightly modified in a way that does not alter
its analytic properties. The resulting covariant identities for the
$\wp_{ijk}$ and $\wp_{ijkl}$ (Appendix 1) differ in detail from
those obtained by Baker (Appendix 2) but are generic and are derived
in a straight forward, economical way with minimal use of computer
algebra and in an algorithmic manner. Because the equivalence of the
two sets of equations is by no means self evident we also specify in
Appendix 2 the transformation between the two definitions of the
$\wp_{ij}$. Given Baker's equations one could have written down the
covariant genus three equations by deducing this simple
transformation by comparing the classical and covariant polar forms.
But this would not have been a test of the machinery nor would it
have provided us with the neat expression for the quadratic, genus
three identities.

By ``minimal use of computer algebra" we mean that the derivation of
the highest weight identities was carried out by hand. A computer
algebra programme was used to implement the actions of $\bf e$ and
$\bf f$ on these highest weight identities in order to check
covariance and to generate the full sets of identities. The other
use of computer algebra, as remarked at the time, was in rearranging
by linear superposition, the identities linear in the $\wp_{ijk}$ in
the genus three case, into the form (\ref{fourterm}).

It may be remarked that Baker's equations are a little simpler than
the covariant ones. From the current point of view this is a
simplification bought at the expense of the more abstract
simplification which incorporates the representation theory. The
drawback of the simplicity is that each identity has to be obtained
independently. The advantage of the marginally more involved
covariant set is that the fifteen identities for the $\wp_{ijkl}$
decompose into sets of nine, five and one elements from each of
which one need only find a single identity using the singularity
analysis, the others following by application of the \emph{raising}
and \emph{lowering operators}, ${\bf e}$ and ${\bf f}$.

Further, the representation theory lays bare a pattern of bones with
further intriguing symmetries that beg further study, particularly
in view of the $\sigma$ function and hyperelliptic addition laws.
However, the most pressing issue now is to apply these methods to
the more difficult non-hyperelliptic curves of low genus.

\section{Appendix 1}
Here we summarize the four-index relations for the covariant
$\wp$-function in the genus three case. $\Delta$ denotes a quadratic
in two index functions:
$\Delta=\wp_{11}\wp_{33}-\wp_{12}\wp_{23}-\wp_{13}^2+\wp_{13}\wp_{22}.$
\begin{eqnarray}
-\wp_{3333}+6\wp_{33}^2&=&8a_6\wp_{33}-8a_7\wp_{23}+a_8(3\wp_{22}-4\wp_{13})\nonumber\\
&&+10(a_4a_8-4a_5a_7+3a_6^2)\nonumber\\
-\wp_{2333}+6\wp_{23}\wp_{33}&=&12a_5\wp_{33}-10a_6\wp_{23}+4a_7(\wp_{22}-3\wp_{13})\nonumber\\
&&+2a_8\wp_{12}\nonumber\\
&&+10(a_3a_8-3a_4a_7+2a_5a_6)\nonumber\\
-\wp_{2233}+4\wp_{23}^2+2\wp_{22}\wp_{33}&=&18a_4\wp_{33}-12a_5\wp_{23}+2a_6(3\wp_{22}-14\wp_{13})\nonumber\\
&&+4a_7\wp_{12}+2a_8\wp_{11}\nonumber\\
&&+8(a_2a_8-a_3a_7-5a_4a_6+5a_5^2)\nonumber\\
-\wp_{2223}+6\wp_{22}\wp_{23}&=&28a_3\wp_{33}-16a_4\wp_{23}+4a_5(3\wp_{22}-14\wp_{13})\nonumber\\
&&+12a_7\wp_{11}\nonumber\\
&&+4(a_1a_8+5a_2a_7-21a_3a_6+15a_4a_5)\nonumber\\
-\wp_{2222}+6\wp_{22}^2-12\Delta&=&48a_3\wp_{33}-32a_3\wp_{23}+32a_4(\wp_{22}-3\wp_{13})\nonumber\\
&&-32a_5\wp_{12}+48a_6\wp_{11}\nonumber\\
&&a_0a_8+24a_1a_7-4a_2a_6-216a_3a_5+195a_4^2\nonumber\\
-\wp_{1333}+6\wp_{13}\wp_{33}&=&3a_4\wp_{33}-10a_6\wp_{13}+4a_7\wp_{12}-a_8\wp_{11}\nonumber\\
&&+3a_2a_8-8a_3a_7+5a_4a_6\nonumber\\
-\wp_{1233}+4\wp_{13}\wp_{23}+2\wp_{12}\wp_{33}&=&4a_3\wp_{33}+2a_4\wp_{23}-20a_5\wp_{13}+6a_6\wp_{12}\nonumber\\
&&+2(a_1a_8-6a_3a_6+5a_4a_5)\nonumber\\
-\wp_{1223}+4\wp_{12}\wp_{23}+2\wp_{13}\wp_{22}&=&6a_2\wp_{33}+4a_3\wp_{23}+2a_4(\wp_{22}-18\wp_{13})\nonumber\\
+2\Delta&&+4a_5\wp_{12}+6a_6\wp_{11}\nonumber\\
&&+\frac12(a_0a_8+16a_1a_7-36a_2a_6-16a_3a_5+35a_4^2)\nonumber\\
-\wp_{1222}+6\wp_{12}\wp_{22}&=&12a_1\wp_{33}+4a_3(3\wp_{22}-14\wp_{13})\nonumber\\
&&-16a_4\wp_{12}+28a_5\wp_{11}\nonumber\\
&&+4(a_0a_7+5a_1a_6-21a_2a_5+15a_3a_4)\nonumber\\
-\wp_{1133}+4\wp_{13}^2+2\wp_{11}\wp_{33}-2\Delta&=&4a_3\wp_{23}-a_4(\wp_{22}-12\wp_{13})+4a_5\wp_{12}\nonumber\\
&&+\frac12(a_0a_8-16a_3a_5+15a_4^2)\nonumber\\
-\wp_{1123}+4\wp_{12}\wp_{13}+2\wp_{11}\wp_{23}&=&6a_2\wp_{23}-20a_3\wp_{13}+2a_4\wp_{12}+4a_5\wp_{11}\nonumber\\
&&+2(a_0a_7-6a_2a_5+5a_3a_4)\nonumber\\
-\wp_{1122}+4\wp_{12}^2+2\wp_{11}\wp_{22}&=&2a_0\wp_{33}+4a_1\wp_{23}+2a_2(3\wp_{22}-14\wp_{13})\nonumber\\
&&-12a_3\wp_{12}+18a_4\wp_{11}\nonumber\\
&&+8(a_0a_6-a_1a_5-5a_2a_4+5a_3^2)\nonumber\\
-\wp_{1113}+6\wp_{11}\wp_{13}&=&-a_0\wp_{33}+4a_1\wp_{23}-10a_2\wp_{13}+3a_4\wp_{11}\nonumber\\
&&+3a_0a_6-8a_1a_5+5a_2a_4\nonumber\\
-\wp_{1112}+6\wp_{11}\wp_{12}&=&2a_0\wp_{23}+4a_1(\wp_{22}-3\wp_{13})\nonumber\\
&&-10a_2\wp_{12}+12a_3\wp_{11}\nonumber\\
&&+10(a_0a_5-3a_1a_4+2a_2a_3)\nonumber\\
-\wp_{1111}+6\wp_{11}^2&=&a_0(3\wp_{22}-4\wp_{13})-8a_1\wp_{12}+8a_2\wp_{11}\nonumber\\
&&+10(a_0a_4-4a_1a_3+3a_2^2)\nonumber
\end{eqnarray}

\section{Appendix 2}
Here we reproduce the genus three equations from Baker's paper
\cite{B1903}. We will denote by $\wp^{\mathfrak B}$ the traditional
genus three $\wp$-function. For ease of comparison we have also
rewritten the $\lambda_i$ coefficients of the monomials in $x$ in
the octic curve in Baker's paper in terms of the $a_i$ used above.

\begin{eqnarray}
{{\wp^{\mathfrak B}}}_{3333}-6{\wp^{\mathfrak B}}_{33}^2&=&28a_6{\wp^{\mathfrak B}}_{33}+8a_7{\wp^{\mathfrak B}}_{23}\nonumber\\
&&+a_8(4{\wp^{\mathfrak B}}_{13}-3{\wp^{\mathfrak B}}_{22})\nonumber\\
&&-35a_4a_8+56a_5a_7\nonumber\\
{\wp^{\mathfrak B}}_{2333}-6{\wp^{\mathfrak B}}_{23}{\wp^{\mathfrak B}}_{33}&=&28a_6{\wp^{\mathfrak B}}_{23}+4a_7(3{\wp^{\mathfrak B}}_{13}-{\wp^{\mathfrak B}}_{22})\nonumber\\
&&+2a_8{\wp^{\mathfrak B}}_{12}-14a_3a_8\nonumber\\
{\wp^{\mathfrak B}}_{2233}-4{\wp^{\mathfrak B}}_{23}^2-2{\wp^{\mathfrak B}}_{22}{\wp^{\mathfrak B}}_{33}&=&28a_5{\wp^{\mathfrak B}}_{23}+28a_6{\wp^{\mathfrak B}}_{13}-4a_7{\wp^{\mathfrak B}}_{12}\nonumber\\
&&-2a_8{\wp^{\mathfrak B}}_{11}-14a_2a_8\nonumber\\
{\wp^{\mathfrak B}}_{2223}-6{\wp^{\mathfrak B}}_{22}{\wp^{\mathfrak B}}_{23}&=&-28a_3{\wp^{\mathfrak B}}_{33}+7-a_4{\wp^{\mathfrak B}}_{23}+56a_5{\wp^{\mathfrak B}}_{13}\nonumber\\
&&-12a_7{\wp^{\mathfrak B}}_{11}\nonumber\\
&&-4a_1a_8-56a_2a_7\nonumber\\
{\wp^{\mathfrak B}}_{2222}-6{\wp^{\mathfrak B}}_{22}^2-12\Delta&=&-84a_2{\wp^{\mathfrak B}}_{33}+56a_3{\wp^{\mathfrak B}}_{23}+70a_4{\wp^{\mathfrak B}}_{22}\nonumber\\
&&+56a_5{\wp^{\mathfrak B}}_{12}-84a_6{\wp^{\mathfrak B}}_{11}\nonumber\\
&&-392a_2a_6+392a_3a_5\nonumber\\
{\wp^{\mathfrak B}}_{1333}-6{\wp^{\mathfrak B}}_{13}{\wp^{\mathfrak B}}_{33}&=&28a_6{\wp^{\mathfrak B}}_{13}-4a_7{\wp^{\mathfrak B}}_{14}+a_8{\wp^{\mathfrak B}}_{11}\nonumber\\
{\wp^{\mathfrak B}}_{1233}-4{\wp^{\mathfrak B}}_{13}{\wp^{\mathfrak B}}_{23}-2{\wp^{\mathfrak B}}_{12}{\wp^{\mathfrak B}}_{33}&=&28a_5{\wp^{\mathfrak B}}_{13}-2a_1a_8\nonumber\\
{\wp^{\mathfrak B}}_{1223}-4{\wp^{\mathfrak B}}_{12}{\wp^{\mathfrak B}}_{23}-2{\wp^{\mathfrak B}}_{13}{\wp^{\mathfrak B}}_{22}&=&70a_4{\wp^{\mathfrak B}}_{13}-8a_1a_7-\frac12a_0a_8\nonumber\\
+2\Delta&&\nonumber\\
{\wp^{\mathfrak B}}_{1222}-6{\wp^{\mathfrak B}}_{12}{\wp^{\mathfrak B}}_{22}&=&-12a_1{\wp^{\mathfrak B}}_{33}+56a_3{\wp^{\mathfrak B}}_{13}+70a_4{\wp^{\mathfrak B}}_{12}\nonumber\\
&&-28a_5{\wp^{\mathfrak B}}_{11}\nonumber\\
&&-112a_1a_6-4a_0a_7\nonumber\\
{\wp^{\mathfrak B}}_{1133}-4{\wp^{\mathfrak B}}_{13}^2-2{\wp^{\mathfrak B}}_{11}{\wp^{\mathfrak B}}_{33}&=&-\frac12a_1a_8\nonumber\\
-2\Delta&&\nonumber\\
{\wp^{\mathfrak B}}_{1123}-4{\wp^{\mathfrak B}}_{12}{\wp^{\mathfrak B}}_{13}-2{\wp^{\mathfrak B}}_{11}{\wp^{\mathfrak B}}_{23}&=&28a_3{\wp^{\mathfrak B}}_{13}-2a_0a_7\nonumber\\
{\wp^{\mathfrak B}}_{1122}-4{\wp^{\mathfrak B}}_{12}^2-2{\wp^{\mathfrak B}}_{11}{\wp^{\mathfrak B}}_{22}&=&-2a_0{\wp^{\mathfrak B}}_{33}-4a_1{\wp^{\mathfrak B}}_{23}+28a_2{\wp^{\mathfrak B}}_{13}\nonumber\\
&&+28a_3{\wp^{\mathfrak B}}_{12}-14a_0a_6\nonumber\\
{\wp^{\mathfrak B}}_{1113}-6{\wp^{\mathfrak B}}_{11}{\wp^{\mathfrak B}}_{13}&=&a_0{\wp^{\mathfrak B}}_{33}-4a_1{\wp^{\mathfrak B}}_{23}+28a_2{\wp^{\mathfrak B}}_{13}\nonumber\\
{\wp^{\mathfrak B}}_{1112}-6{\wp^{\mathfrak B}}_{11}{\wp^{\mathfrak B}}_{12}&=&-2a_0{\wp^{\mathfrak B}}_{23}+4a_1(3{\wp^{\mathfrak B}}_{13}-{\wp^{\mathfrak B}}_{22})\nonumber\\
&&+28a_2{\wp^{\mathfrak B}}_{12}-14a_0a_5\nonumber\\
{\wp^{\mathfrak B}}_{1111}-6{\wp^{\mathfrak B}}_{11}^2&=&a_0(4\wp_{13}-3{\wp^{\mathfrak B}}_{22})+8a_1{\wp^{\mathfrak B}}_{12}\nonumber\\
&&+28a_2{\wp^{\mathfrak B}}_{11}\nonumber\\
&&-35a_0a_4+56a_1a_3\nonumber
\end{eqnarray}

The Baker equivalent of the $5\times 5$ matrix $h$ we will call
$h^{\mathfrak B}$ and since our covariant form is to be replaced by
the classical polar form $h^{\mathfrak B}$ will be given by:

\begin{tiny}
\begin{equation}
\left[\begin{array}{ccccc}
a_0&4a_1&-2{\wp^{\mathfrak B}}_{11}&-2{\wp^{\mathfrak B}}_{12}&-2{\wp^{\mathfrak B}}_{13}\\
4a_1&28a_2+4{\wp^{\mathfrak B}}_{11}&28a_3+2{\wp^{\mathfrak B}}_{12}&-2{\wp^{\mathfrak B}}_{22}+4{\wp^{\mathfrak B}}_{13}&-2{\wp^{\mathfrak B}}_{23}\\
-2{\wp^{\mathfrak B}}_{11}&28a_3+2{\wp^{\mathfrak B}}_{12}&70a_4+4{\wp^{\mathfrak B}}_{22}-4{\wp^{\mathfrak B}}_{13}&28a_5+2{\wp^{\mathfrak B}}_{23}&-2{\wp^{\mathfrak B}}_{33}\\
-2{\wp^{\mathfrak B}}_{12}&-2{\wp^{\mathfrak B}}_{22}+4{\wp^{\mathfrak B}}_{13}&28a_5+2{\wp^{\mathfrak B}}_{23}&28a_6+4{\wp^{\mathfrak B}}_{33}&4a_7\\
-2{\wp^{\mathfrak B}}_{13}&-2{\wp^{\mathfrak B}}_{23}&-2{\mathfrak
B}_{33}&4a_7&a_8
\end{array}
\right]
\end{equation}
\end{tiny}

Consequently
\begin{eqnarray}
{\wp^{\mathfrak B}}_{11}&=&\wp_{11}-3a_2\nonumber\\
{\wp^{\mathfrak B}}_{12}&=&\wp_{12}-2a_3\nonumber\\
{\wp^{\mathfrak B}}_{13}&=&\wp_{13}-\frac12a_4\nonumber\\
{\wp^{\mathfrak B}}_{22}&=&\wp_{22}-9a_4\nonumber\\
{\wp^{\mathfrak B}}_{23}&=&\wp_{23}-2a_5\nonumber\\
{\wp^{\mathfrak B}}_{33}&=&\wp_{33}-3a_6\nonumber\\
\end{eqnarray}

Substitution for either $\wp$ or $\wp^{\mathfrak B}$ does indeed
transform the two sets of equations into one another.

\section{Acknowledgements}

I should like to extend thanks to Victor Enolskii and Chris Eilbeck,
long term collaborators in this area, for their continued interest,
optimism and prodding and to Jan Sanders for clearing up some
distracting issues. Thanks are also very much due to the referees
for their perceptive remarks and suggestions for improvements.

\bibliographystyle{amsplain}

\end{document}